# General Solution to Sequential Linear Conformable Fractional Differential Equations With Constant Coefficients

Emrah Ünal[a], Ahmet Gökdoğan[b], Ercan Çelik[c]

[a] Department of Elementary Mathematics Education, Artvin Çoruh University, 08100 Artvin, Turkey

emrah.unal@artvin.edu.tr

[b] Department of Mathematical Engineering, Gümüşhane University, 29100 Gümüşhane, Turkey,

gokdogan@gumushane.edu.tr

[c] Department of Mathematics, Atatürk University, 25100 Erzurum, Turkey,

ecelik@atauni.edu.tr

**Abstract**

In this work, we give the general solution sequential linear conformable fractional differential equations in the case of constant coefficients for $\alpha \in (0,1]$. In homogeneous case, we use a fractional exponential function which generalizes the corresponding ordinary function. In non-homogeneous case, we present to fractional the method of variation of parameters for a particular solution of sequential linear conformable fractional differential equations.

**MSC:** 26A33

**Key Words:** Sequential Linear Conformable Fractional Differential Equations, Conformable Fractional Derivative, Fractional Method of Variation of Parameters

## 1. Introduction

Though L'Hospital has suggested the concept of fractional derivative long back (17th century), still its reflection is being found in several researches of recent centuries. For defining the fractional derivative, most researchers preferred to apply the integral form. Caputo and Riemann-Liouville are the two frequently used definitions. For understanding these definitions in detail, readers are referred to read [1-3].

Recently fractional derivative and associated integral have been freshly defined by Khalil and colleagues [4]. New definitions take the advantages from the limit form as used in the regular derivatives. This new theory has been improved by Abdeljawad [5]. For instance, he gives definitions of left and right conformable fractional derivatives, Taylor power series representation and Laplace transformation of certain functions, fractional integration by parts formulas, chain rule and Gronwall inequality.

In short time, a lot of studies about new fractional derivative definition have been presented. Some works in this field are with regard to the power series solutions around an ordinary point and a regular-singular point homogeneous sequential linear conformable fractional differential equations of order $2\alpha$ in the case of variable coefficients [6,7,8,9], conformable fractional fourier series [10], boundary value problems for conformable fractional differential equations [11,12] and existence and uniqueness theorems for sequential linear conformable fractional differential equations [13].

While the general solution has been found for Riemann-Liouville fractional derivative, $\alpha$-exponential function which is Mittag-Leffler-type function is used [14]. However we present the general solution using a fractional exponential function for sequential linear conformable fractional differential equations in the case of constant coefficients.

The Method of Variation of Parameters are not applicable due to a noticeable lack of basic properties in the Riemann-Liouville derivative. However this case for conformable fractional derivative is not valid. Finally, In this work, we present the fractional method of variation parameters to derive a particular solution.

## 2. Conformable Fractional Calculus

**Definition 2.1.** $f:[0,\infty) \to \mathbb{R}$ let a function. Then for all $t > 0$, the conformable fractional derivative of $f$ of order $\alpha$ is defined as

$$T_\alpha(f)(t) = \lim_{\varepsilon \to 0} \frac{f(t+\varepsilon t^{1-\alpha})-f(t)}{\varepsilon}$$

where $\alpha \in (0,1)$.

If $f$ is $\alpha$-differentiable in some $(0,a), a > 0$ and $\lim_{t \to 0^+} f^{(\alpha)}(t)$ exists, then define

$$f^{(\alpha)}(0) = \lim_{t \to 0^+} f^{(\alpha)}(t).$$

**Theorem 2.1.** Let $\alpha \in (0,1]$ and $f,g$ be $\alpha$-differentiable at a point $t > 0$. Then

(1) $T_\alpha(af + bg) = aT_\alpha(f) + bT_\alpha(g)$, for all $a,b \in \mathbb{R}$.

(2) $T_\alpha(t^p) = pt^{p-\alpha}$ for all $p \in \mathbb{R}$.

(3) $T_\alpha(\lambda) = 0$ for all constant functions $f(t) = \lambda$.

(4) $T_\alpha(fg) = T_\alpha(f)g + fT_\alpha(g)$.

(5) $T_\alpha(f/g) = \frac{T_\alpha(f)g - fT_\alpha(g)}{g^2}$.

(6) In addition, If $f$ is differentiable, then $T_\alpha(f(t)) = t^{1-\alpha}\frac{df}{dt}$.

Additionaly, conformable fractional derivatives of certain functions as follow:

(i) $T_\alpha\left(\sin\frac{1}{\alpha}t^\alpha\right) = \cos\frac{1}{\alpha}t^\alpha$.
(ii) $T_\alpha\left(\cos\frac{1}{\alpha}t^\alpha\right) = -\sin\frac{1}{\alpha}t^\alpha$.
(iii) $T_\alpha\left(e^{\frac{1}{\alpha}t^\alpha}\right) = e^{\frac{1}{\alpha}t^\alpha}$.

**Definition 2.2.** $f:[a,\infty) \to \mathbb{R}$ let a function. Then for all $t > a$, $\alpha \in (0,1)$ the "conformable fractional integral" of $f$ of order $\alpha$ is

$$(I_\alpha^a f)(t) = \int_a^t f(x) d_\alpha(x) = \int_a^t x^{\alpha-1} f(x) dx$$

where the integral is the usual Riemann improper integral.

**Theorem 2.2.** Let $f:[a,\infty) \to R$ be any continuous function and $0 < \alpha \leq 1$. Then for all $t > a$

$$T_\alpha I_\alpha^a f(t) = f(t).$$

**Theorem 2.3.** Let $f, g:[a,b] \to \mathbb{R}$ be two functions such that $f, g$ is differentiable. Then

$$\int_a^b f(t) T_\alpha(g)(t) d_\alpha(t) = fg|_a^b - \int_a^b g(t) T_\alpha(f)(t) d_\alpha(t).$$

## 3. General Solution in the Homogeneous Case

In this section we introduce a method, analogous to that for the ordinary case. Let $y$ be $n$ times $\alpha$-differentiable function for $\alpha \in (0,1]$. The most general sequential linear homogeneous conformable fractional differential equation with constant coefficients is

$$^n T_\alpha y + p_{n-1}{}^{n-1} T_\alpha y + \cdots + p_2\, ^2 T_\alpha y + p_1 T_\alpha y + p_0 y = 0. \tag{1}$$

where $^n T_\alpha y = T_\alpha T_\alpha \ldots T_\alpha y$, $n$ times, and the coefficients $p_0, p_1, \ldots, p_{n-1}$ are real constants. Left-hand of equation (1) rewrite by

$$L_\alpha[y] = \left(^n T_\alpha + p_{n-1}{}^{n-1} T_\alpha + \cdots + p_2\, ^2 T_\alpha + p_1 T_\alpha + p_0 \right) y. \tag{2}$$

If the equation (1) have $y_1(t), y_2(t), \ldots, y_n(t)$ that are linearly independent solutions, general solution is

$$y = c_1 y_1(t) + c_2 y_2(t) + \cdots + c_n y_n(t).$$

where $c_1, c_2, \ldots, c_n$ are arbitrary constants [13].

**Lemma 3.1.** Let $L_\alpha[.]$ is a linear operator with constant coefficients and $\alpha \in (0,1]$. For $t > 0$

$$L_\alpha\left[e^{\frac{r}{\alpha} t^\alpha}\right] = P_n(r) e^{\frac{r}{\alpha} t^\alpha}$$

where $r$ is real or complex constant and $P_n(r) = r^n + p_{n-1} r^{n-1} + \cdots + p_0$.

**Proof :** Conformable derivatives of $y = e^{\frac{r}{\alpha} t^\alpha}$ are

$$T_\alpha y = r e^{\frac{r}{\alpha} t^\alpha}, \quad ^2 T_\alpha y = r^2 e^{\frac{r}{\alpha} t^\alpha}, \ldots, \quad ^n T_\alpha y = r^n e^{\frac{r}{\alpha} t^\alpha}. \tag{3}$$

If $y = e^{\frac{r}{\alpha}t^\alpha}$ and the equations (3) substitute in $L_\alpha[y] = (\,^nT_\alpha + p_{n-1}\,^{n-1}T_\alpha + \cdots + p_2\,^2T_\alpha + p_1 T_\alpha + p_0)y$, then

$$L_\alpha\left[e^{\frac{r}{\alpha}t^\alpha}\right] = (r^n + p_{n-1}r^{n-1} + \cdots + p_0)e^{\frac{r}{\alpha}t^\alpha},$$

$$L_\alpha\left[e^{\frac{r}{\alpha}t^\alpha}\right] = P_n(r)e^{\frac{r}{\alpha}t^\alpha}$$

is obtained. Hence, the proof is completed.

As for the usual case, we shall seek the solution of the equation (1) in the form $y = e^{\frac{r}{\alpha}t^\alpha}$, where $r$ is a real or complex constant. It follows from the equation (2) and Lemma 3.1 that

$$L_\alpha\left[e^{\frac{r}{\alpha}t^\alpha}\right] = P_n(r)e^{\frac{r}{\alpha}t^\alpha} = 0.$$

$P_n(r) = r^n + p_{n-1}r^{n-1} + \cdots + p_0$ is called as the characteristic polynomial. For all $r$ we have $e^{\frac{r}{\alpha}t^\alpha} \neq 0$. Hence, $P_n(r) = r^n + p_{n-1}r^{n-1} + \cdots + p_0 = 0$ is obtained. Here,

$$r^n + p_{n-1}r^{n-1} + \cdots + p_0 = 0 \qquad (4)$$

is called as the characteristic equation.

**Lemma 3.2:** If $r$ is a root of the characteristic equation (4), then

$$\frac{\partial}{\partial r}\left[L_\alpha\left[e^{\frac{r}{\alpha}t^\alpha}\right]\right] = L_\alpha\left[\frac{\partial}{\partial r}e^{\frac{r}{\alpha}t^\alpha}\right]$$

and

$$\frac{\partial^l}{\partial r^l}e^{\frac{r}{\alpha}t^\alpha} = \left(\frac{t^\alpha}{\alpha}\right)^l e^{\frac{r}{\alpha}t^\alpha}.$$

**Proof:** $L_\alpha[\,.\,]$ is linear as is seen from Theorem 4.3 in [13]. Also, $\frac{\partial}{\partial r}$ is linear as is known the classical derivative. Hence, we can written

$$\frac{\partial}{\partial r}\left[L_\alpha\left[e^{\frac{r}{\alpha}t^\alpha}\right]\right] = L_\alpha\left[\frac{\partial}{\partial r}e^{\frac{r}{\alpha}t^\alpha}\right].$$

Moreover, by the help of classical derivative, it seen that

$$\frac{\partial^l}{\partial r^l}e^{\frac{r}{\alpha}t^\alpha} = \frac{t^{l\alpha}}{\alpha^l}e^{\frac{r}{\alpha}t^\alpha}.$$

**Lemma 3.3:** Let $r_1$ is a root which multiplicity $\mu_1$ of the characteristic equation (4). Then, for $l = 0, 1, \dots \mu_1 - 1$, the functions

$$y_{1,l}(t) = \left(\frac{t^\alpha}{\alpha}\right)^l e^{\frac{r_1}{\alpha}t^\alpha}$$

are solutions of equation (1).

**Proof:** By using Lemma 3.2, we can written

$$\left\{L_\alpha\left[\frac{\partial^l}{\partial r^l}e^{\frac{r}{\alpha}t^\alpha}\right]\right\}_{r=r_1} = \left\{\frac{\partial^l}{\partial r^l}\left[L_\alpha\left[e^{\frac{r}{\alpha}t^\alpha}\right]\right]\right\}_{r=r_1}.$$

Considering the equation $L_\alpha\left[e^{\frac{r}{\alpha}t^\alpha}\right] = P_n(r)e^{\frac{r}{\alpha}t^\alpha}$ and using classical Leibniz rule,

$$\left\{L_\alpha\left[\frac{\partial^l}{\partial r^l}e^{\frac{r}{\alpha}t^\alpha}\right]\right\}_{r=r_1} = \sum_{j=0}^{l}\binom{l}{j}\left[\frac{\partial^{l-j}}{\partial r^{l-j}}\left(e^{\frac{r}{\alpha}t^\alpha}\right)\right]_{r=r_1}\frac{\partial^j}{\partial r^j}[P_n(r)]_{r=r_1}$$

is obtained. Since $r_1$ is a root of multiplicity $\mu_1$ of the characteristic equation (4), for $j = 0, 1, \dots, \mu_1 - 1$, we have

$$\left[\frac{\partial^j}{\partial r^j}P_n(r)\right]_{r=r_1} = 0.$$

From Lemma 3.2, $\frac{\partial^l}{\partial r^l}e^{\frac{r}{\alpha}t^\alpha} = \left(\frac{t^\alpha}{\alpha}\right)^l e^{\frac{r_1}{\alpha}t^\alpha} = y_{1,l}(t)$. Hence,

$$L_\alpha[y_{1,l}(t)] = 0$$

is obtained.

**Corollary 3.1:** Let $r_1, r_2, \dots, r_k$ are distinct roots of, respectively, multiplicity $\mu_1, \mu_2, \dots, \mu_k$ of the characteristic equation (4). So, linearly independent set of solutions to equation (1) for these roots is as following:

$$\bigcup_{m=1}^{k}\left\{\left(\frac{t^\alpha}{\alpha}\right)^l e^{\frac{r_m}{\alpha}t^\alpha}\right\}_{l=0}^{\mu_m-1}.$$

**Proof:** Corollary 4.1 is seen by Lemma 3.3 and Theorem 4.5 in [13].

**Lemma 3.4:** Let $r_1$ ve $\overline{r_1}$ ($r_1 = \theta + i\beta, \beta \neq 0$) are complex roots which multiplicity $\sigma_1$ of the characteristic equation (4). Hence, for $l = 0, 1, \ldots \sigma_1 - 1$, the functions

$$y_{1,l}(t) = \left(\frac{t^\alpha}{\alpha}\right)^l e^{\frac{\theta}{\alpha}t^\alpha}\left[\cos\left(\frac{\beta}{\alpha}t^\alpha\right) + i\sin\left(\frac{\beta}{\alpha}t^\alpha\right)\right]$$

and

$$y_{2,l}(t) = \left(\frac{t^\alpha}{\alpha}\right)^l e^{\frac{\theta}{\alpha}t^\alpha}\left[\cos\left(\frac{\beta}{\alpha}t^\alpha\right) - i\sin\left(\frac{\beta}{\alpha}t^\alpha\right)\right]$$

are linearly independent solutions of equation (1).

**Proof:** Since $r_1 = \theta + i\beta$ is a root of multiplicity $\sigma_1$ of the characteristic equation (4), according Lemma 3.3, the functions

$$y_{1,l}(t) = \left(\frac{t^\alpha}{\alpha}\right)^l e^{\frac{\theta+i\beta}{\alpha}t^\alpha}$$

are solutions of the equation (1). Analogously, for $\overline{r_1} = \theta - i\beta$, the functions

$$y_{2,l}(t) = \left(\frac{t^\alpha}{\alpha}\right)^l e^{\frac{\theta-i\beta}{\alpha}t^\alpha}$$

are also solutions of the equation (1). For these solutions, if the Euler's identity:

$$e^{i\frac{\beta}{\alpha}t^\alpha} = \cos\frac{\beta}{\alpha}t^\alpha + i\sin\frac{\beta}{\alpha}t^\alpha$$
$$e^{-i\frac{\beta}{\alpha}t^\alpha} = \cos\frac{\beta}{\alpha}t^\alpha - i\sin\frac{\beta}{\alpha}t^\alpha$$

is used, the solutions given in Lemma 3.4 are get.

**Corollary 3.2:** Let $\{r_m, \overline{r_m}\}_{m=1}^p$, $r_m = \theta_m + i\beta_m, \beta_m \neq 0$, are distinct $2p$ roots of multiplicity $\{\sigma_m\}_{m=1}^p$ of the characteristic equation (4). In this case, for these roots, elements of following sets are linearly independent solutions of the equation (1):

$$\bigcup_{m=1}^{p}\left\{\left(\frac{t^\alpha}{\alpha}\right)^l e^{\frac{\theta_m}{\alpha}t^\alpha}\left[\cos\left(\frac{\beta_m}{\alpha}t^\alpha\right)+i\sin\left(\frac{\beta_m}{\alpha}t^\alpha\right)\right]\right\}_{l=0}^{\sigma_m-1}$$

and

$$\bigcup_{m=1}^{p}\left\{\left(\frac{t^\alpha}{\alpha}\right)^l e^{\frac{\theta_m}{\alpha}t^\alpha}\left[\cos\left(\frac{\beta_m}{\alpha}t^\alpha\right)-i\sin\left(\frac{\beta_m}{\alpha}t^\alpha\right)\right]\right\}_{l=0}^{\sigma_m-1}.$$

These sets have linearly independent $2\sum_{m=1}^{p}\sigma_m$ solutions.

**Proof:** Proof is seen by Lemma 3.4 and Theorem 4.5 in [13].

**Theorem 3.1:** Let $\{r_j\}_{j=1}^{k}$ are distinct $k$ roots of multiplicity $\{\mu_j\}_{j=1}^{k}$ and $\{\lambda_j, \bar{\lambda}_j\}_{j=1}^{p}$, $\lambda_j = \theta_j + i\beta_j, \beta_j \neq 0$, are distinct $2p$ roots which multiplicity $\{\sigma_j\}_{j=1}^{p}$ of the equation (1) such as $\sum_{j=1}^{k}\mu_j + 2\sum_{j=1}^{p}\sigma_j = n$. So, linearly independent set of solutions to the equation (1) is union of following sets:

$$\bigcup_{m=1}^{k}\left\{\left(\frac{t^\alpha}{\alpha}\right)^l e^{\frac{r_m}{\alpha}t^\alpha}\right\}_{l=0}^{\mu_m-1},$$

$$\bigcup_{m=1}^{p}\left\{\left(\frac{t^\alpha}{\alpha}\right)^l e^{\frac{\theta_m}{\alpha}t^\alpha}\left[\cos\left(\frac{\beta_m}{\alpha}t^\alpha\right)+i\sin\left(\frac{\beta_m}{\alpha}t^\alpha\right)\right]\right\}_{l=0}^{\sigma_m-1}$$

and

$$\bigcup_{m=1}^{p}\left\{\left(\frac{t^\alpha}{\alpha}\right)^l e^{\frac{\theta_m}{\alpha}t^\alpha}\left[\cos\left(\frac{\beta_m}{\alpha}t^\alpha\right)-i\sin\left(\frac{\beta_m}{\alpha}t^\alpha\right)\right]\right\}_{l=0}^{\sigma_m-1}.$$

**Proof:** The proof is seen by Corollary 3.1, Corollary 3.2 and Theorem 4.5 in [13].

**Example 3.1.** $^2T_\alpha y + 4T_\alpha y + 3y = 0.$

The characteristic equation of the above equation is

$$r^2 + 4r + 3 = 0.$$

The roots are

$$r = -3 \text{ and } r = -1.$$

Hence, the general solution is

$$y(t) = c_1 e^{\frac{-3}{\alpha}t^\alpha} + c_2 e^{\frac{-1}{\alpha}t^\alpha}.$$

**Example 3.2.** $^2T_\alpha y - 10 T_\alpha y + 25 y = 0.$

The characteristic equation of the above equation is

$$r^2 - 10r + 25 = (r-5)^2 = 0.$$

The roots are $r_{1,2} = 5$. Hence, the general solution is

$$y = \left(c_1 + c_2 t^\alpha\right) e^{\frac{5}{\alpha}t^\alpha}.$$

**Example 3.3.** $^2T_\alpha y + T_\alpha y + y = 0.$

The characteristic equation of the above equation is

$$r^2 + r + 1 = 0.$$

The roots are

$$r_1 = -\frac{1}{2} + i\frac{\sqrt{3}}{2} \text{ and } r_2 = -\frac{1}{2} - i\frac{\sqrt{3}}{2}.$$

Hence, the general solution is

$$y(t) = e^{-\frac{1}{2\alpha}t^\alpha}\left(c_1 \cos\frac{\sqrt{3}}{2\alpha}t^\alpha + c_2 \sin\frac{\sqrt{3}}{2\alpha}t^\alpha\right).$$

## 4. Method of Variation of Parameters for Conformable Fractional Calculus

In this section, we apply method of variation of parameters to derivate the particular solution of equation

$$^n T_\alpha y + p_{n-1}{}^{n-1}T_\alpha y + \cdots + p_2\,^2T_\alpha y + p_1 T_\alpha y + p_0 y = q(t) \tag{5}$$

where $y$ is $n$ times $\alpha$-differentiable function for $\alpha \in (0,1]$ at a point $t > 0$.

**Theorem 4.1.** Let $u(t)$ be a function which is solution of homogenous case of equation (5) and given by

$$u(t) = \sum_{i=1}^n c_i y_i(t). \tag{6}$$

Then particular solution of the equation (5) is

$$v(t) = \sum_{i=1}^n c_i(t) y_i(t)$$

where $c_1(t), c_2(t), \ldots, c_n(t)$ provide following systems of equations

$$\sum_{i=1}^{n} c_i^{(\alpha)}(t) y_i(t) = 0$$

$$\sum_{i=1}^{n} c_i^{(\alpha)}(t) y_i^{(\alpha)}(t) = 0$$

$$\vdots$$

$$\sum_{i=1}^{n} c_i^{(\alpha)}(t)^{n-2} T_\alpha\, y_i(t) = 0$$

$$\sum_{i=1}^{n} c_i^{(\alpha)}(t)^{n-1} T_\alpha\, y_i(t) = q(t).$$

**Proof:** Now we shall seek the solution of the equation (5) in the form

$$v(t) = \sum_{i=1}^{n} c_i(t) y_i(t)$$

If we calculate conformable derivative of $v(t)$ for $\alpha \in (0,1]$, then we get

$$T_\alpha v(t) = \sum_{i=1}^{n} c_i(t) y_i^{(\alpha)}(t) + \sum_{i=1}^{n} c_i^{(\alpha)}(t) y_i(t).$$

Let first condition is $\sum_{i=1}^{n} c_i^{(\alpha)}(t) y_i(t) = 0$. In this case, we obtain

$$T_\alpha v(t) = \sum_{i=1}^{n} c_i(t) y_i^{(\alpha)}(t).$$

If we calculate conformable derivative of $T_\alpha v(t)$ for $\alpha \in (0,1]$, then we get

$$^2T_\alpha v(t) = \sum_{i=1}^{n} c_i(t)\, ^2T_\alpha y_i(t) + \sum_{i=1}^{n} c_i^{(\alpha)}(t) y_i^{(\alpha)}(t).$$

Let second condition is $\sum_{i=1}^{n} c_i^{(\alpha)}(t) y_i^{(\alpha)}(t) = 0$. So, we have

$$^2T_\alpha v(t) = \sum_{i=1}^{n} c_i(t)\, ^2T_\alpha y_i(t).$$

By continuing in this way;

$$^{n-1}T_\alpha v(t) = \sum_{i=1}^{n} c_i(t)^{n-1} T_\alpha y_i(t) + \sum_{i=1}^{n} c_i^{(\alpha)}(t)^{n-2} T_\alpha y_i(t)$$

is obtained. Let $(n-1)^{th}$ condition is $\sum_{i=1}^{n} c_i^{(\alpha)}(t)^{n-2} T_\alpha y_i(t) = 0$. Hence, it is obtained that

$$^{n-1}T_\alpha v(t) = \sum_{i=1}^{n} c_i(t)^{n-1} T_\alpha y_i(t).$$

And finally, calculating conformable fractional derivative of above equation for $\alpha \in (0,1]$, we can write the following equation

$$^nT_\alpha v(t) = \sum_{i=1}^n c_i(t)\,^nT_\alpha y_i(t) + \sum_{i=1}^n c_i^{(\alpha)}(t)\,^{n-1}T_\alpha y_i(t).$$

If we substitute $v(t), T_\alpha v(t),\ ^2T_\alpha v(t)\ ...,\ ^nT_\alpha v(t)$ in the equation (5), then we have

$$\sum_{i=1}^n c_i^{(\alpha)}(t)^{n-1}T_\alpha y_i(t) + \sum_{i=1}^n c_i(t)\left[p_o y_i(t) + p_1 y_i^{(\alpha)}(t) + \cdots + \,^nT_\alpha y_i(t)\right] = q(t).$$

Because $y_1(t), y_2(t), \ldots, y_n(t)$ are solutions of equation (6), we get

$$\sum_{i=1}^n c_i(t)\left[p_o y_i(t) + p_1 y_i^{(\alpha)}(t) + \cdots + \,^{(n)}y_i^{(\alpha)}(t)\right] = 0.$$

Hence,

$$\sum_{i=1}^n c_i^{(\alpha)}(t)^{n-1}T_\alpha y_i(t) = q(t)$$

is obtained. The last equation is $n^{th}$ condition. The system of equations formed by these conditions is as following:

$$\sum_{i=1}^n c_i^{(\alpha)}(t) y_i(t) = 0$$

$$\sum_{i=1}^n c_i^{(\alpha)}(t) y_i^{(\alpha)}(t) = 0$$

$$\vdots$$

$$\sum_{i=1}^n c_i^{(\alpha)}(t)^{n-2}T_\alpha\, y_i(t) = 0$$

$$\sum_{i=1}^n c_i^{(\alpha)}(t)^{n-1}T_\alpha\, y_i(t) = q(t).$$

If $c_1(t), c_2(t), \ldots, c_n(t)$ provide system of equations which are built by conditions, then function $v(t) = \sum_{i=1}^n c_i(t) y_i(t)$ is a particular solution of the equation (5).

**Example 4.1.** $^2T_\alpha y + 4T_\alpha y + 3y = q(t)$.

a) Let $q(t) = e^{2t^\alpha}$. For $v(t) = c_1(t) e^{\frac{-3}{\alpha}t^\alpha} + c_2(t) e^{\frac{-1}{\alpha}t^\alpha}$, the system of equations which are built by conditions is

$$c_1^{(\alpha)}(t) e^{\frac{-3}{\alpha}t^\alpha} + c_2^{(\alpha)}(t) e^{\frac{-1}{\alpha}t^\alpha} = 0$$

$$-3 c_1^{(\alpha)}(t) e^{\frac{-3}{\alpha}t^\alpha} - c_2^{(\alpha)}(t) e^{\frac{-1}{\alpha}t^\alpha} = e^{2t^\alpha}.$$

Solving the above system of equations, we get $c_1^{(\alpha)}(t) = -\frac{1}{2}e^{\frac{2\alpha+3}{\alpha}t^\alpha}$, $c_2^{(\alpha)}(t) = \frac{1}{2}e^{\frac{2\alpha+1}{\alpha}t^\alpha}$. For this values, if we calculate conformable integral of this values for $\alpha \in (0,1]$, then we have
$c_1(t) = -\frac{1}{4\alpha+6}e^{\frac{2\alpha+3}{\alpha}t^\alpha}$, $c_2(t) = \frac{1}{4\alpha+2}e^{\frac{2\alpha+1}{\alpha}t^\alpha}$. The particular solution $v(t)$ is

$$v(t) = -\frac{1}{4\alpha+6}e^{\frac{2\alpha+3}{\alpha}t^\alpha} \cdot e^{\frac{-3}{\alpha}t^\alpha} + \frac{1}{4\alpha+2}e^{\frac{2\alpha+1}{\alpha}t^\alpha} \cdot e^{\frac{-1}{\alpha}t^\alpha}$$

$$= \frac{1}{4\alpha^2+8\alpha+3}e^{2t^\alpha}.$$

b) Let $q(t) = 2t^{2\alpha} + t^\alpha - 3$. In this case, the system of equations is

$$c_1^{(\alpha)}(t)e^{\frac{-3}{\alpha}t^\alpha} + c_2^{(\alpha)}(t)e^{\frac{-1}{\alpha}t^\alpha} = 0$$

$$-3c_1^{(\alpha)}(t)e^{\frac{-3}{\alpha}t^\alpha} - c_2^{(\alpha)}(t)e^{\frac{-1}{\alpha}t^\alpha} = 2t^{2\alpha} + t^\alpha - 3.$$

If we solve this system of equations, we obtain

$$c_1(t) = -\frac{1}{3}t^{2\alpha}e^{\frac{3}{\alpha}t^\alpha} + \frac{4\alpha-3}{18}t^\alpha e^{\frac{3}{\alpha}t^\alpha} + \frac{-4\alpha^2+3\alpha+27}{54}e^{\frac{3}{\alpha}t^\alpha},$$

$$c_2(t) = t^{2\alpha}e^{\frac{1}{\alpha}t^\alpha} + \frac{1-4\alpha}{2}t^\alpha e^{\frac{1}{\alpha}t^\alpha} + \frac{4\alpha^2-2\alpha-6}{2}e^{\frac{1}{\alpha}t^\alpha}.$$

Hence, we obtain particular solution $v(t)$ as following:

$$v(t) = \frac{2}{3}t^{2\alpha} + \frac{3-16\alpha}{9}t^\alpha + \frac{52\alpha^2-12\alpha-27}{27}.$$

c) Let $q(t) = \sin 2t^\alpha$. In this case, the system of equations is

$$c_1^{(\alpha)}(t)e^{\frac{-3}{\alpha}t^\alpha} + c_2^{(\alpha)}(t)e^{\frac{-1}{\alpha}t^\alpha} = 0$$

$$-3c_1^{(\alpha)}(t)e^{\frac{-3}{\alpha}t^\alpha} - c_2^{(\alpha)}(t)e^{\frac{-1}{\alpha}t^\alpha} = \sin 2t^\alpha.$$

From the system of equations, we find $c_1^{(\alpha)}(t) = -\frac{1}{2}(\sin 2t^\alpha)e^{\frac{3}{\alpha}t^\alpha}$, $c_2^{(\alpha)}(t) = \frac{1}{2}(\sin 2t^\alpha)e^{\frac{1}{\alpha}t^\alpha}$. Using integration by parts for conformable fractional derivative, we have

$$c_1(t) = -\frac{3}{8\alpha^2+18}\sin 2t^\alpha e^{\frac{3}{\alpha}t^\alpha} + \frac{\alpha}{4\alpha^2+9}\cos 2t^\alpha e^{\frac{3}{\alpha}t^\alpha},$$

$$c_2(t) = \frac{1}{8\alpha^2+2}\sin 2t^\alpha e^{\frac{1}{\alpha}t^\alpha} - \frac{\alpha}{4\alpha^2+1}\cos 2t^\alpha e^{\frac{1}{\alpha}t^\alpha}.$$

So, we get

$$v(t) = \frac{-8\alpha}{16\alpha^2+40\alpha+9}\cos 2t^\alpha + \frac{-4\alpha^2+3}{16\alpha^4+40\alpha^2+9}\sin 2t^\alpha.$$

d) Let $q(t) = e^{2t^\alpha}t^\alpha$. We can write the system of equations as following:

$$c_1^{(\alpha)}(t)e^{\frac{-3}{\alpha}t^\alpha} + c_2^{(\alpha)}(t)e^{\frac{-1}{\alpha}t^\alpha} = 0$$

$$-3c_1^{(\alpha)}(t)e^{\frac{-3}{\alpha}t^\alpha} - c_2^{(\alpha)}(t)e^{\frac{-1}{\alpha}t^\alpha} = e^{2t^\alpha}t^\alpha.$$

Solving the system of equations, we have

$$c_1(t) = -\frac{1}{4\alpha+6}e^{\frac{2\alpha+3}{\alpha}t^\alpha}t^\alpha + \frac{\alpha}{2(2\alpha+3)^2}e^{\frac{2\alpha+3}{\alpha}t^\alpha},$$

$$c_2(t) = \frac{1}{4\alpha+2}e^{\frac{2\alpha+1}{\alpha}t^\alpha}t^\alpha + \frac{\alpha}{2(2\alpha+1)^2}e^{\frac{2\alpha+1}{\alpha}t^\alpha}.$$

Then, the particular solution

$$v(t) = \frac{e^{2t^\alpha}t^\alpha}{4\alpha^2+8\alpha+3} - \frac{4\alpha^2+4\alpha}{(4\alpha^2+8\alpha+3)^2}e^{2t^\alpha}$$

is obtained.

e) Let $q(t) = e^{-4t^\alpha}$. For $\alpha \neq 3/4$ and $\alpha \neq 1/4$, the system of equations is

$$c_1^{(\alpha)}(t)e^{\frac{-3}{\alpha}t^\alpha} + c_2^{(\alpha)}(t)e^{\frac{-1}{\alpha}t^\alpha} = 0$$

$$-3c_1^{(\alpha)}(t)e^{\frac{-3}{\alpha}t^\alpha} - c_2^{(\alpha)}(t)e^{\frac{-1}{\alpha}t^\alpha} = e^{-4t^\alpha}.$$

From the system of equations, we have $c_1(t) = \frac{1}{8\alpha-3}e^{\frac{3-4\alpha}{\alpha}t^\alpha}$, $c_2(t) = \frac{1}{2-8\alpha}e^{\frac{1-4\alpha}{\alpha}t^\alpha}$. In this case, we obtain the particular solution as following:

$$v(t) = \frac{e^{-4t^\alpha}}{16\alpha^2-16\alpha+3}.$$

For $\alpha = 3/4$ and $v(t) = c_1(t)e^{-4t^{3/4}} + c_2(t)e^{-\frac{4}{3}t^{3/4}}$, the system of equations is obtained as

$$c_1^{(3/4)}(t)e^{-4t^{3/4}} + c_2^{(3/4)}(t)e^{-\frac{4}{3}t^{3/4}} = 0$$

$$-3c_1^{(3/4)}(t)e^{-4t^{3/4}} - c_2^{(3/4)}(t)e^{-\frac{4}{3}t^{3/4}} = e^{-4t^{3/4}}.$$

From this system of equations, we have $c_1(t) = -\frac{2}{3}t^{3/4}$, $c_2(t) = -\frac{1}{4}e^{-\frac{8}{3}t^{3/4}}$. The particular solution

$$v(t) = -\frac{2}{3}t^{3/4}e^{-4t^{3/4}} - \frac{1}{4}e^{-4t^{3/4}}$$

is obtained.

Similarly, for $\alpha = 1/4$ and $v(t) = c_1(t)e^{-12t^{1/4}} + c_2(t)e^{-4t^{1/4}}$, the system of equations is

$$c_1^{(1/4)}(t)e^{-12t^{1/4}} + c_2^{(1/4)}(t)e^{-4t^{1/4}} = 0$$

$$-3c_1^{(1/4)}(t)e^{-12t^{1/4}} - c_2^{(1/4)}(t)e^{-4t^{1/4}} = e^{-4t^{1/4}}.$$

From this system of equations, we find $c_1(t) = -\frac{1}{4}e^{8t^{1/4}}$, $c_2(t) = 2t^{1/4}$. The particular solution

$$v(t) = 2t^{1/4}e^{-4t^{1/4}} - \frac{1}{4}e^{-4t^{1/4}}$$

is obtained.

## 5. Conclusion

In this work, using fractional exponential function, we give the general solution to sequential linear homogeneous and non-homogeneous differential equation of conformable fractional of order with constant coefficients for $0 < \alpha \leq 1$. In contrast to the Riemann-Liouville fractional derivatives, it has been found that the method of variation of parameters can be applied for conformable fractional derivatives. Additionally, it has appeared that the results obtained in this work correspond to results which are obtained in ordinary cases.